\documentclass[a4paper,reqno,11pt]{amsart}
\usepackage{amsmath}\usepackage{amssymb}\usepackage{latexsym}
\usepackage{graphicx}
\usepackage{longtable}\usepackage{array}

\font\black=msbm10 scaled\magstep1

\def\field #1{\hbox{{\black #1}}}

\def\pd#1#2{\frac{\partial#1}{\partial#2}}
\def\set#1{\{\,#1\,\}}                 
\def\>#1{{\bf #1}}                      

\def\eq#1{{\begin{equation} #1 \end{equation}}}

\def\mitad{\frac{1}{2}}

\def\a{\alpha}
\def\R{{\hbox{{\field R}}}}                 
\def\C{{\hbox{{\field C}}}}              
\def\N{{\hbox{{\field N}}}}              

\def\im{\hbox{{\rm Im}}}                          

\def\rango{\hbox{{\rm rank}}}
\def\rank{\hbox{{\rm rank}}}

\newtheorem{theorem}{Theorem}

\newcommand{\nc}{\newcommand}                        
\nc{\rnc}{\renewcommand}

\rnc{\a}{\alpha} \rnc{\b}{\beta} \rnc{\d}{\delta} \nc{\D}{\Delta} \nc{\e}{\varepsilon}
\nc{\g}{\gamma} \nc{\G}{\Gamma} \rnc{\l}{\lambda} \rnc{\L}{\Lambda}  \nc{\n}{\nabla}
\nc{\var}{\varphi} \nc{\ro}{\rho} \nc{\s}{\sigma} \nc{\Sig}{\Sigma} \rnc{\t}{\tau} \nc{\te}{\theta}
\rnc{\o}{\omega} \rnc{\O}{\Omega} \nc{\z}{\zeta}
\rnc{\P}{\Phi} \nc{\Te}{\Theta}
\def\eqy#1{{\begin{eqnarray} #1 \end{eqnarray}}}
\def\eqi#1{{\begin{eqnarray*} #1 \end{eqnarray*}}}
\def\matrix#1#2{\left[\begin{array}{#1} #2 \end{array}\right]}

\begin{document}

\title{A Numerical Algorithm for Singular Optimal LQ  Control Systems}

\author{Marina Delgado--T\'ellez}\address{Marina Delgado--T\'ellez: Departamento de Matem\'aticas\\Universidad
Carlos III de Madrid\\
Avda. de la Universidad 30, Legan\'{e}s 28911\\Madrid\\Spain}\email{mdelgado@math.uc3m.es}

\author{Alberto Ibort}\address{Alberto Ibort: Departamento de Matem\'aticas\\Universidad Carlos III de Madrid\\
Avda. de la Universidad 30, Legan\'{e}s 28911\\Madrid\\Spain}\email{albertoi@math.uc3m.es}

\thanks{Research partially supported by MEC grant MTM2004-07090-C03-03. SIMUMAT-CM, CP05/CP06-CM-UCIIIM}


\begin{abstract}
A numerical algorithm to obtain the consistent conditions satisfied by singular arcs for singular
linear-quadratic optimal control problems is presented.  The algorithm  is based on the 
presymplectic constraint algorithm (PCA) by Gotay-Nester \cite{Go78,Vo99} that allows to solve presymplectic hamiltonian systems and that provides a geometrical framework to the Dirac-Bergmann theory of constraints for singular Lagrangian systems \cite{Di49}.   The numerical implementation of the algorithm is based on the singular value decomposition that, on each step allows to construct a semi-explicit system.  Several examples and experiments are discussed, among them a family of arbitrary large singular LQ systems with index 3 and a family of examples of arbitrary large index, all of them exhibiting stable behaviour.
\end{abstract}
\maketitle
\tableofcontents

\newpage
\section{Introduction}

Singular problems in optimal control problems and in the calculus of variations have been widely
considered from different perspectives.  A singular perturbation approach has been often the preferred
approach to singular optimal control theory (see for instance \cite{Ko84,Ko86,Om78} and references
therein). More recently Jurdjevic has proposed a different viewpoint by introducing Lie-theoretic
methods in the problem \cite{Ju97}.   A similar problem in the calculus of variations has had a very
different history mainly due to the physical insight derived from its appearance in mechanics and
field theory, (see for instance Cari\~nena \cite{Ca90}). 
P.A.M. Dirac took a brilliant approach introducing a
recursive analysis to extract the integrable or solvable part of the system \cite{Di49}. 
Dirac's approach consists
essentially on a recursive consistency scheme where the original system of implicit Euler-lagrange equations obtained from the extremal conditions are restricted to a smaller and smaller subset on state
space by imposing the existence of at least one solution to the problem passing through them.  Such approach has been adapted to the problem of optimal control by Volkaert \cite{Vo}, L\'opez and Mart\'{\i}nez \cite{Lo00}, Guerra \cite{Gu01} and Delgado and Ibort \cite{De03a} by transforming it into a descriptor system.    A geometric version of Dirac's constraint algorithm was presented by Gotay and Nester \cite{Go78} for presymplectic systems.   This algorithm was extended and generalized later on to more general situations (see for instance \cite{Me05}) arriving to essentially the same geometrical algorithm devised by Rabier and Rehinbold \cite{Ra94} to deal with Differential-Algebraic equations (DAE's).  A variety of ideas and techniques have been developed along the years to analyze DAE's (see for instance \cite{Br89}, \cite{Ra02} or the recent book \cite{Ku06}) and
recently various of these approaches have been applied to singular optimal control problems (see \cite{Ku97}, \cite{Ba06} for instance). 
However the PCA algorithm was proposed for the first time in control theory by  \cite{Vo99}, and L\'opez and Mart\'{\i}nez \cite{Lo00},  Cort\'es,  de Leon,  Mart\'{\i}n de Diego and  Mart\'{\i}nez \cite{Co03} and it has been discussed by Delgado and Ibort \cite{De03a,De03b,De04} in the specific context of singular control theory. Another form of this algorithm was
applied to the LQ singular case by  Guerra \cite{Gu01}.   

In this paper we will concentrate on singular linear-quadratic systems.  We will transform such
algorithm into a linear algebra problem that can be treated numerically.
In \cite{Ku97} the methods developed by Kunkel and Merhmann based on the construction of normal forms were applied to predictor systems and as a particular instance, to singular LQ optimal problems More recently \cite{Ba06} have extended these ideas incorporating to the analysis LQ systems with index 3.    Campbell's differential arrays method can also be applied to singular optimal control problems providing both a description of consistent initial conditions and  the differential equation giving the solutions to the problem \cite{Br89}.    A full account of the analysis of linear differential-algebraic equations with variable coefficients is presented that could be applied to general singular LQ systems in \cite{Ku06}.  We have chosen however to implement the PCA algorithm to solve singular LQ optimal control systems because, apart that it can lead to computational improvements with respect to the general procedures above, it preserves the structure of the system.   The linear DAE's that are obtained from the analysis of singular LQ systems carry a presymplectic structure (inherited from the canonical symplectic structure on the space of states and coestates).     Such structure emerges at the end on the set of consistent states and coestates inducing there a (pre)symplectic structure.  The reduced equations are hamiltonian with respect to such structure.    The algorithm that will be implemented in the present paper will compute the consistent initial states and the corresponding consistent coestates preserving the structure of the problem even though we will leave the analysis of the reduced hamiltonian equations to a continuation of this work.    

The numerical implementation of the algorithm is carefully discussed and a family of examples and experiments are analyzed.   We will discuss experiments with large matrices and lower index (2 and 3) and experiments with large index $n-1$, where $n$ is the dimension of the state space exhibiting all of them a remarkable stable behaviour of the algorithm.   We must point it out however that in the present state of the analysis of the numerical algorithm it is not possible to prove its backwards stability because it contains repeated products of matrices spoiling such possibility.   More refined versions of it will be analyzed elsewhere, with a different handling of rank conditions that will improve its efficiency and that will allow for a rigorous error analysis.

The paper will be organized as follows.  First we will set up in Section \ref{sing-lq} the basic notions for singular optimal problems and the particular instance of LQ systems that we are going to discuss.   Section \ref{kron} will be devoted to review the PCA algorithm from the slightly wider perspective of quasilinear implicit differential equations and we will discuss the relation of the recursive index of the algorithm to the standard index.  In Section \ref{pcasect}  we will describe the linear algebraic
algorithm corresponding to this problem and in Section \ref{experiments}  we will describe its
numerical stability properties by means of various experiments.

\section{Constraint algorithms for singular LQ systems}\label{sing-lq}

As it was stated in the introduction we will concentrate on the study of LQ optimal control systems. 
We will discuss the problem of finding $C^1$-piecewise smooth curves $\gamma(t)=(x(t),u(t))$
satisfying the linear control equation:
\begin{equation}\label{control}
    \dot{x}^i=A^i_j\,x^j+B^i_a\,u^a,
\end{equation}and minimizing the objective functional:
\begin{equation}
    \label{cost} S(\gamma)=\int_{t_0}^TL(x(t),u(t))\,d t,
\end{equation}
where the quadratic Lagrangian $L$ has the form:
\begin{equation}
    L(x,u)=\frac{1}{2}Q_{ij}\,x^i\,x^j+N_{ia}\,x^i\,u^a+\frac{1}{2}R_{ab}\,u^a\,u^b,
\end{equation}
subjected to fixed endpoints conditions: $x(t_0)=x_0$, $x(T)=x_T$ (however we must point it out that the chosen endpoints conditions are not going to be relevant for the analysis to follow, hence they can be replaced by more general ones without altering the results presented in this paper).

The coordinates $x^i$, $i=1,\ldots,n$, describe points $x\in \R^n$ in state space and $u^a$,
$a=1,\ldots,m$, are control coordinates  defined on the linear control space $\R^m$. The matrices
$A$, $B$, $Q$, $N$ and $R$ will be considered to be constant for simplicity.  Again we must stress that the algorithm that we are going to discuss can be applied without difficulty to the time-dependent situtation, however the numerical implementation becomes much more involved, so we have chosen to discuss it just for time-independent systems in order to gain clarity and a better understanding of the experiments and their stability properties.

It is well known that normal etxtremals  to this problem are provided by Pontryagin's Maximum
Principle \cite{Po62}:

The curve $\gamma(t)=(x(t),\,u(t))$ is a normal extremal trajectory if there exists a lifting
$(x(t), p(t))$ of $x(t)$ to the costate space $\R^n\times\R^n$ satisfying Hamilton's
equations:
\begin{equation}\label{pmp}
    \dot{x}^i=\frac{\partial H}{\partial p_i}, \quad \quad
    \dot{p}_i=-\frac{\partial H}{\partial x^i},
\end{equation}  
where $H$ is  Pontryagin's Hamiltonian function:
\begin{equation}\label{ham}
    H(x,p,u) = p_i\,(A^i_j\,x^j+B^i_a\,u^a)-L(x,u),
\end{equation}
and the set of conditions:
\begin{equation}\label{liga}
    \phi^{(1)}_{\ \ a}(x,p,u):= \frac{\partial H}{\partial u^a}
    =p_i\,B^i_a-N_{ia}\,x^i-R_{ab}\,u^b=0.
\end{equation}

We will call conditions in Eq. (\ref{liga}) primary
constraints.  Thus, trajectories solution
to the optimal control problem must lie in the linear submanifold:
\begin{equation}\label{m1}
    M_1=\{(x,p,u)\in M_0|~\phi^{(1)}_{\ \ a}(x,p,u) = 0 \},
\end{equation}where $M_0 = \set{(x,p,u)\in\R^{2n+m}}$ denotes the total  space of the system.

If $(x(t),p(t),u(t))$ is a solution of the optimal problem, then its derivative must satisfy:
\eqy{\label{derx}\dot{x}^i&=&\frac{\partial H}{\partial p_i}(x,p,u)= A^i_j\,x^j+B^i_a\,u^a,\\
     \label{derp}\dot{p}_i&=&-\frac{\partial H}{\partial x^i}(x,p,u)= -p_jA^j_i+Q_{ij}x^j+N_{ia}u^a,\\
     \label{deru}\dot{u}^a&=&C^a(x,p,u),}
together with $\dot{\phi}^{(1)}_{\ \ a}(x,p,u)=0$,  this is,
$$\pd{\phi^{(1)}_{\ \ a}}{x^i}\pd{H}{p_i}-\pd{\phi^{(1)}_{\ \ a}} {p_i}\pd{H}{x^i}+
     \pd{\phi^{(1)}_{\ \ a}}{u^b}\,C^b=0.$$

A simple computation shows us that if the system is regular, that is, if the matrix:
$$R_{ab}=\frac{\partial^2 H} {\partial u^a\partial u^b}$$
is invertible in $M_1$,  then there
exists an optimal feedback condition solving Equation (\ref{liga}),
\eq{\label{feedback} u^b=(R^{-1})^{ab}\,(p_iB^i_{\,a}-N_{ia}x^i).}%
Then we obtain for Eq. (\ref{deru}):
\eqi{\dot{u}^b&=&(R^{-1})^{ab}(\dot{p}_iB^i_a-N_{ia}\dot{x}^i)=\\
&=&(R^{-1})^{ab}\left((-p_jA^j_i+Q_{ij}x^j+N_{id}u^d)B^i_a-N_{ia}(A^i_j\,x^j+B^i_d\,u^d)\right).}
Notice that in this case $\dot{\phi}^{(1)}_{\ \ a}$ vanishes automatically on $M_1$.

However, for singular optimal LQ  systems, this is, when $R_{ab}$ is not an invertible matrix, it
may occur that at some points $(x_0, p_0,u_0)$ satisfying the primary constraint Eq. (\ref{liga}), that  solutions of (\ref{pmp}) starting at them will not be contained in $M_1$ for $t> 0$ for any $u$.   Because Eqs. (\ref{pmp})-(\ref{liga}) must be satisfied along optimal paths, we must consider only as initial conditions only those points for which there is at least a solution of (\ref{pmp})
contained in $M_1$ starting from them.   Such subset is defined by the following conditions:%
\begin{equation}\label{recurr}
\mathrm{There~exists~} C \mathrm{~such~that~} \phi^{(2)}_{\ \ a}:=\dot{\phi}^{(1)} _{\ \ a}=0,
\end{equation}%
where the derivative is taken in the direction of Eq. (\ref{pmp}) and $\dot{u}=C$. The subset obtained, that we will denote by  $M_2$, is again a linear submanifold of $M_1$ (this is also true in the time-dependent case) and we
shall denote the functions defining  $M_2$ in $M_1$ by $\phi^{(2)} _{\ \ a}$.   We will call them
secondary constraints.    Notice that, in the case of more general endpoint conditions, that would involve considering end-time conditions on the coestate variables $p_i$, the same condition Eq. (\ref{recurr}) would apply by changing now the derivative of $p_i$ by the derivative along Eq. (\ref{pmp}) obtained by time reversing $t \mapsto -t$.

Clearly the argument goes on and we will obtain in this form a family of linear
submanifolds defined recursively as follows:%
\begin{equation} M_{k+1}=\{(x,p,u)\in M_k \mid \, \exists\ C~{\mbox{such that}}~Ê
\phi^{(k+1)} _{\ \ \ \ \ a}:=\dot{\phi}^{(k)} _{\ \ a}\,(x,p,u) =0\}. \end{equation}
Eventually the recursion will stop and $M_r = M_{r+1} = M_{r+2} = \dots$, for certain finite
$r$. We will call the number of steps $r$ of
the algorithm before it stabilizes the recursive index of the problem.  In this way we obtain an invariant linear submanifold,%
\begin{equation}M_{\infty}=\bigcap_{k\geq 0}\,M_k,
\end{equation}%
that will be called the final constraint submanifold of the problem and  by construction it consists on the set of consistent initial condition for the DAE Eqs. (\ref{pmp})-(\ref{liga}).    The geometrical analysis of this algorithm shows that this number $r$ does not depends on the coordinate system and constitutes a intrinsic property of the system.   Notice also that for general
singular systems, this index may vary in principle from point to point.   However this is not the case for
singular linear systems as we will see in next section where we will discuss for completeness the relation between the recursive index $r$ and the index of linear DAEs.

As it was stressed in the introduction, this algorithm constitutes both an adapted version of the reduction algorithm for DAEs \cite{Ra94,Ra02}  and the  \emph{Presymplectic Constraint Algorithm} (PCA).     The DAE system above however has an additional structure because it is a presymplectic system.    The PCA algorithm not only determines the consistent initial conditions for the corresponding DAE, but in addition it provides the explicit form of the reduced Hamiltonian equations by computing the so called Dirac brackets of the system.  Such brackets are obtained directly from the sequence of $k$-ary constraints 
$\phi^{(k)}_{\ \ a}$.   In this sense this algorithm is structured and preserves the main structure of the system along its steps.    We will just proceed to the numerical implementation of this algorithm in its present form leaving the construction (and integration) of the reduced dynamical system to subsequent articles.

\section{The Kronecker and recursive index for linear DAEs}\label{kron}

As it was indicated before we will show for completeness the relation of the recursive index with the strangeness and differential indexes of the general theory of DAE's.    Because of the simplicity of the problem at hand it will suffice to compute the Kroneker index of the matrix pencil defining the system.
Thus, we consider an (autonomous) quasilinear differential--algebraic equation (DAE) of the form
\eqy{\label{daegen}A(x)\dot{x}=B(x),}
$x\in \R^n = M_0$, where $A(x),\,B(x)\colon\R^n\to F$ and $F$ is an auxiliary
linear space (see \cite{Gr04} and the references therein for a geometrical treatment of DAEs).
Because the DAE above has no additional structures the PCA becomes particularly simply and completely equivalent to the differential reduction by Rabier and Reinholt \cite{Ra96}.   Even more, it can be written in the extremely simple form that follows.
If $A(x)$ is regular, we can solve explicitly $\dot{x}$ and the DAE becomes an ordinary differential equation.  If $A(x)$ is not regular for some $x\in M_0$ we have
to impose the constraints algorithm. We shall define $M_1$ as the set of points in $M_0$ such that 
$B(x)\in\im\,A(x)$.  Hence if $x\in M_1$  then
$\forall\mu\in\ker\,A(x)^*\subset F^*$ it must be satisfied $\left<\mu,B(x)\right>=0$, where  $A^*$
denotes the adjoint application of $A$. In general, \eq{\label{reck1}M_{k+1}=\{x\in
M_k|~B_k(x)\in\im\,A_k(x)\},}where $B_k=\left.B\right|_{M_k};$ $A_k=\left.A\right|_{M_k}.$ We
obtain again \eq{\label{mk1}M_{k+1}=\{x\in
M_k|~\left<\mu,B_k(x)\right>=0,~\forall\mu\in\ker\,A_k(x)^*\}.}

In the particular case of constant linear systems, Equation (\ref{daegen}) becomes
\eqy{\label{dae}A\cdot\dot{x}=B\cdot x,}
where $x,~\dot{x}\in\R^n;~A,~B\in\R^{n\times n},$ and the
submanifold $M_1$ defined by (\ref{reck1}) is given by
$$M_1=\{x\in\R^n|~B\cdot x\in\im\,A\},$$
but this is equivalent to say that for all $z$ verifying
$z^T\cdot A=0$ then $z^T\cdot B\cdot x=0$. If $z_a,~a=1,\ldots, m,$ is a basis of  $\ker\,(A^T)$,
we  can construct the family of primary constraints as%
\eqi{\phi^{(1)}_{\ \ a}(x)=z_a^T\cdot B\cdot x.}%
If we define the matrix $C^{(1)}:=[z_1|\cdots|z_m]$,  then  the family of linear constraints
$\{\phi^{(1)}_{\ \ a}\}_{a=1}^m$ is equivalent to the following matrix equation
\eq{C^{(1)T}\cdot B\cdot x=0.}%
Now the condition $x\in M_1$ is equivalent to $x\in\ker\,(C^{(1)T}\cdot B)$.
Denoting by $A_k$ as the restriction of $A$ to $M_k,$ we obtain that the linear manifold $M_{k+1}$
is defined recursively as the set of points $x\in M_k$ such that
\eq{C^{(k+1)T}\cdot B\cdot x=0,}
where the columns of $C^{(k+1)}$ 
generates the kernel of the matrix $A_k$.


Let $A\cdot\dot{x}=B\cdot x$ be a constant implicit differential--algebraic system.  We will recall
that the system is regular in the Kronecker sense if the matrix pencil $A\lambda-B$, $\l\in\C$, is
regular, i.e., if the set of solutions of the characteristic equation
$p(\lambda)=\det(A\lambda-B)=0$ is finite.

Moreover  if the pencil is regular in the Kronecker sense, then there exists regular  matrices $E$,
$F$ such that
\cite{Gan59}:\eq{EAF=\matrix{c|c}{I&0\\\hline0&N};~~EBF=\matrix{c|c}{W&0\\\hline0&I},}where $N$ is a
nilpotent matrix  with index $\nu$.  In such case we will say that the index of the
implicit differential-algebraic system  given by Equation (\ref{dae}) is $\nu$ (see \cite{Ku06} for a thourough discussion of the subject).   The relation of the recursive index to the index of the DAE is given by the following result:

\theorem{\label{teoindex} The index $\nu$ of the regular constant implicit
differential-algebraic system $A\dot{x}=Bx$ coincides with the number of steps of the recursive
constraint algorithm minus one: \eq{\mathit{number~of~steps~}=r=\nu+1.} }

\proof
If Equation (\ref{dae}) is a regular system  in the sense of Kronecker, there exist  regular matrices
$E$ and $F$ such that the system becomes \eqy{\dot{y}_1&=&Wy_1,\\\label{ny}N\dot{y}_2&=&y_2,}%
where $N$ is nilpotent with index $\nu$ (i.e. $N^{\nu}\neq 0,~N^{\nu+1}=0$) and $x=Fy$, together
with $y=\matrix{c}{y_1\\\hline y_2}$.

Now we only need to consider the nilpotent part of the system above, Equation (\ref{ny}).
Applying the recursive constraint algorithm to it, we obtain the primary constraint submanifold
$M_1$,  given by the primary constraints
$$\phi^{(1)}=C^{(1)T}y_2,$$where the columns of $C^{(1)}$ generates the kernel of $N$.

Moreover,  Equation  (\ref{ny}) implies that $y_2\in M_1$ if and only if exists $z\in M_0$ such
that $y_2=Nz$, i.e., $y_2\in\im\, N$. If $y_2(t)$ is a curve solution of the equation, it will mean
that $y_2(t)=Nz(t)$, hence  $\dot{y}_2(t)=N\dot{z}(t)$ and then $\dot{y}_2\in\im\, N$.
So we have that $y_2=N\dot{y}_2=NN\dot{z}=N^2\dot{z}.$ Then $y_2\in M_2$ if and only if
$y_2\in\im\, N^2,$  and this will happen if and only if
$$\phi^{(2)}=C^{(2)T}y_2=0,$$where $C^{(2)}$ generates the kernel of $N^2,$ thus
$$M_2=\left\{y_2\in M_1|~C^{(2)T}y_2=0\right\},~\mathrm{Lin}\{\mathrm{Col}(C^{(2)})\}=\ker N^2.$$

If we proceed recursively, we can observe that $y_2\in M_{k+1}$ if and only if $y_k\in\im\, N^k$
and this will take place if and only if $\phi^{(k)}=C^{(k)T}y_2=0,$ where the columns of $C^{k}$
generates the kernel of the matrix $N^k$. As the matrix $N$ is nilpotent with index $\nu$ we have
that
$$\ker\,N\subsetneqq\ker\,N^2\subsetneqq\ker\,N^3\subsetneqq\cdots\subsetneqq\ker\,N^{\nu}
\subsetneqq\ker\,N^{\nu+1}=\R^{m}.$$

So, in each step, the matrix $C^{(k)}$  contains the previous one, $C^{(k-1)}$, as a submatrix,
$$C^{(k)}=[C^{(k-1)}|~*~].$$%
Notice that in the $\nu$--th step we obtain that
$$y_2\in\im\,N^{\nu}\Longrightarrow y_2=N^{\nu}z,$$so,  the next step is
$y_2=N\dot{y}_2=NN^{\nu}\dot{z}=0$ and the final constraint submanifold is
$M_{\nu+1}=M_{\infty}=\{y_2=0\},$ given by the constraints $y_2=0.$
\endproof

\section{A numerical linear algebra algorithm for singular LQ systems}\label{pcasect}

Now we will adapt the general recursive constraint algorithm stated in Section \ref{sing-lq} to the
particular instance of singular LQ control systems.    Notice that the description of the constrainsts algorithm done in Section \ref{sing-lq}  does not involves the use of the presymplectic structure, thus it is a plain constraints algorithm for a DAE.   We will follow this approach here instead of using the full PCA algorithm because here we are just addressing the problem of determining the set of consistent initial conditions for our problem.    The idea that we are going to follow to implement the algorithm is to transform at each step the implicit problem we have into a semi-explicit system.  Then, we will not only get the set of constraints defining the set of consistent initial condition but we will have at the same time the set of explicit reduced equations of the system.    In this form the algorithm will not provide the Hamiltonian structure of the equations though as it was pointed it out before.

We will discuss now the basic idea
of the algorithm from the matrix analysis perspective.
We will write down all coordinates as column vectors: $(x,p)\in\R^n\times\R^n$ and $u\in\R^m$. The
control equation and the lagrangian density have the form already described in Section \ref{sing-lq},
now written in matrix notation reads as
\eqy{\dot{x} & = & Ax+Bu, \\L & = & \mitad x^TQx+x^TNu+\mitad u^TRu,}with $A,~Q\in\R^{n\times n}$,
$B,~N\in\R^{n\times m}$ and $R\in\R^{m\times m}$. The Pontryagin's Hamiltonian becomes:
\eq{H(x,p,u)=p^TAx+p^TBu-\mitad x^TQx-x^TNu-\mitad u^TRu.}%
Using these notations the equations of motion (\ref{derx})-(\ref{deru})  become:
\eqy{ \dot{x} & = & \pd{H}{p^T}=Ax+Bu, \\ \dot{p} & = & -\pd{H}{x^T}=-A^Tp+Qx+Nu,\\ \dot{u} & = &C,}%
and the column primary constraint vector  is given by:
\eq{\phi^{(1)} :=  -N^Tx+B^Tp-Ru.}

We must notice that when applying the recursive constraint algorithm to the problem above, that all the
constraints thus obtained will be linear. Then we can writet them at each step $k$ of the algorithm
as:
\eq{\label{ligak}\phi^{(k)}(x,p,u) := \sigma^{(k)}x + \beta^{(k)}p + \rho^{(k)}u ,}%
with matrices $\sigma^{(k)}, \beta^{(k)} \in \R^{r_k\times n}$, $\rho^{(k)} \in \R^{r_k\times m}$, for some
$r_k\in\N$.  Notice that:
$$ \sigma^{(1)} = -N^T, \quad \beta^{(1)} = B^T, \quad \rho^{(1)} = - R .$$%
The matrix equations
$$\phi^{(1)} = 0 ; \ldots ; \phi^{(k)} = 0,$$%
define the linear manifolds $M_k$ obtained by applying the recursive constraint algorithm. Thus the
matrices $\sigma^{(k)}$, $\beta^{(k)}$, $\rho^{(k)}$ completely characterize the constraints
$\phi^{(k)}$. We will store these matrices in a block structured matrix $\P$ whose $k$--th row,
$\P(k,:)$, will be given by $ \matrix{c|c|c}{ \sigma^{(k)} & \beta^{(k)} & \rho^{(k)} }.$

If  $R$ is a regular  matrix, then there exists an optimal feedback and the control variables are
uniquely determined. However, if $R$ is singular we must apply the recursive constraint algorithm.
As we have discussed above, the first step of the algorithm amounts to study the stability of the
primary constraint $\phi^{(1)}$, i.e., to determine for which points $(x,p,u)$ there exists a
vector $C\in\R^m$ satisfying
\eq{\label{ligalq}\dot{\phi}^{(1)}=\sigma^{(1)}\dot{x}+\beta^{(1)}\dot{p}+\rho^{(1)}C=0,}%
with $C=\dot{u}$.  Because $\rho^{(1)}$ is singular, the linear system obtained from it
$$ \rho^{(1)}C= b(x,p,u) ,$$%
with $b(x,p,u) = -\sigma^{(1)}\dot{x}-\beta^{(1)}\dot{p}$, will not always have solution. However,
given the dependence of $(x, p, u)$ on the inhomogeneous part of the linear equation, we can
determine for which values of them the system will have solution.  For those points $(x,p,u)$ such
that a solution exists, i.e.,  such that there exists $C\in\R^m$ with $\dot{\phi}^{(1)}(x,p,u) = 0$, we will
obtain a partial optimal feedback.  The remaining equations will impose further conditions on the
points $(x,p,u)$ where we can expect to find solutions to the original optimal control problem.
That part will constitute what we were calling before the secondary constraints of the problem.

We will obtain this separation between the partial optimal feedback and the secondary constraints by
using the singular value decomposition (SVD) of $\rho^{(1)}$.  There exists two unique orthogonal
matrices, $U^{(1)},~V^{(1)}$, and real numbers  $s_1^{(1)}\geq s_2^{(1)}\geq\dots\geq s_{r_1}^{(1)}
>0,$  the singular values of $\rho^{(1)}$, such that:
$$\rho^{(1)}=U^{(1)}\,\Sigma^{(1)}\,V^{(1)\,T}=U^{(1)}\, \left[\begin{array} {c c c|c}
s_1^{(1)}& & & \\ &\ddots& &0\\ & &s_{r_1}^{(1)}&\\ \hline &0& &0
\end{array}\right]\,V^{(1)\,T} .$$%
Redefining the variables $u^{(1)}=V^{(1)\,T}u$, then
$\dot{u}^{(1)}=V^{(1)\,T}\dot{u}=V^{(1)\,T}C=C^{(1)},$ and
\eqy{0&=&U^{(1)\,T}\,\dot{\phi}^{(1)}=
U^{(1)\,T}\left(\sigma^{(1)}\dot{x}+\beta^{(1)}\dot{p}\right)+\Sigma^{(1)}C^{(1)}=\nonumber\\
\label{feedliga} &=&U^{(1)\,T}\left(\sigma^{(1)}\dot{x}+\beta^{(1)}\dot{p}\right)+
\left[\begin{array}{c|c} \Sigma^{(1)}\,_{r_1}&0\\\hline0&0\end{array}\right]C^{(1)}=0,}%
where $\Sigma^{(1)}\,_{r_1}=\mathrm{diag}(s_1^{(1)},\ldots,s_{r_1}^{(1)}).$%
We will split $C^{(1)}$ as $[C_{r_1}^{(1)},C_{m-r_1}^{(1)}]^T,$ where $C_{r_1}^{(1)}$ are the
components $1,\ldots,r_1,$ of $C^{(1)}$  and $C_{m-r_1}^{(1)}$  the $r_1,\ldots,m,$ ones. We then
get a partial feedback for $C_{r_1}^{(1)}$ and the new constraint $\phi^{(2)}$:
\eq{\label{feed2} [~I_{r_1}|~~~~0~~~~]\,U^{(1)\,T}\left(\sigma^{(1)}\dot{x}+
\beta^{(1)}\dot{p}\right)+\Sigma_{r_1}^{(1)}C_{r_1}^{(1)}=0}
and,
\eqy{\label{nofeed2}\phi^{(2)}&:=& [~0~~|~I_{m-r_1}~]\,U^{(1)\,T}\left(\sigma^{(1)}\dot{x}+\beta^{(1)}\dot{p}\right)=\\
\nonumber &=&[~0~|~I_{m-r_1}~]\,U^{(1)\,T}\left[(\sigma^{(1)}A+ \beta^{(1)}Q)x+(-\beta^{(1)}A^T)p+
(\sigma^{(1)}B+\beta^{(1)}N)u\right]=\\
\nonumber &=& \sigma^{(2)}x+\beta^{(2)}p+\rho^{(2)}u=0.}%
Iterating the process, given the $k$--th constraint $\phi^{(k)}$ defined in (\ref{ligak}), the SVD
of the matrix $\rho^{(k)}$ provides orthogonal matrices $U^{(k)}$, $V^{(k)}$, and the new constraint:
\eqi{\phi^{(k+1)}:=\sigma^{(k+1)}x+\beta^{(k+1)}p+\rho^{(k+1)}u,}
with,
\eqy{\label{rec1}\sigma^{(k+1)}&=&[~0~|~I_{m-r_k}~]\,U^{(k)\,T}\left(\sigma^{(k)}A+
\beta^{(k)}Q\right),\\
\label{rec2}\beta^{(k+1)}&=&[~0~|~I_{m-r_k}~]\,U^{(k)\,T}\left(-\beta^{(k)}A^T \right),\\
\label{rec3}\rho^{(k+1)} &=&[~0~|~I_{m-r_k}~]\,U^{(k)\,T}\left(\sigma^{(k)}B+\beta^{(k)}N\right).}%
The recursive relations above can be solved explicitly.     First we shall denote as
$U^k:=[~0~|~I_{m-rk}~]\,U^{(k)\,T}$ for each $k$, and then consider the simplified recurrence relations:
\eqy{\label{sigmaktilda} \widetilde{\sigma}^{(k+1)} & = & \widetilde{\sigma}^{(k)}A + \widetilde{\beta}^{(k)}Q\\
\label{betaktilda} \widetilde{\beta}^{(k+1)} & = & - \widetilde{\beta}^{(k)}A^T \\
\label{rhoktilda} \widetilde{\rho}^{(k+1)} & = &
\widetilde{\sigma}^{(k)}B + \widetilde{\beta}^{(k)}N}
Before solving this set of conditions, let us compute the relation of $\widetilde{\sigma}$, $\widetilde{\beta}$ and $\widetilde{\rho}$ with $\sigma$, $\beta$ and $\rho$ respectively.
First we will compute such relation for $\beta$. Notice that:
\eq{\beta^{(1) }= \widetilde{\beta}^{(1)}, \quad \beta^{(2)} = U^1\left(-\beta^{(1)}A^T\right)=U^1\widetilde{\beta}^{(2)} ,}
then for $k\geq 2$:
\eqy{\beta^{(k)}&=&U^k\left(-\beta^{(k-1)}A^T\right) =  \nonumber \\ &=&U^{k-1}\cdots
U^2U^1\left(-\widetilde{\beta}^{(k-1)}A^T\right)= U^{k-1}\cdots U^2U^1\widetilde{\beta}^{(k)}.}
A similar computation shows that: 
\eq{\sigma^{(k+1)}=U^k\cdots U^2U^1\widetilde{\sigma}^{(k+1)}, \quad k\geq 1.}
Finally, when computing $\rho^{(k+1)}$ we obtain the same result.  
\eq{\rho^{(k+1)} = U^k\cdots
U^2U^1\left(\widetilde{\sigma}^{(k)}B+\widetilde{\beta}^{(k)}N \right) = U^k\cdots
U^2U^1\widetilde{\rho}^{(k+1)} .}
Now expanding Eqs. (\ref{sigmaktilda})-(\ref{rhoktilda}) we obtain the explicit expression for the matrices $\widetilde{\beta}^{(k)}$:
\eq{ \widetilde{\beta}^{(1)} = B^T, \quad \widetilde{\beta}^{(k+1)} = (-1)^{k}B^T(A^T)^k, \quad k\geq 1 .}
For the matrices $\widetilde{\sigma}^{(k)}$ we get:
\eq{ \widetilde{\sigma}^{(1)} = -N^T, \quad \widetilde{\sigma}^{(k+1)} = -N^TA^k + B^T \left[\sum_{i=0}^{k-1}(-1)^i(A^T)^i Q A^{k-1-i}\right] ,}
and, finally for the matrices $\widetilde{\rho}^{(k)}$ we obtain:
\eqy{ \widetilde{\rho}^{(1)} ~Ê~Ê~ & = & -R, \quad  \widetilde{\rho}^{(2)} = -N^T B + B^T N, \nonumber \\     \widetilde{\rho}^{(k+1)} & = &  -N^TA^{k-1}B
+(-1)^{k-1}B^T(A^T)^{k-1}N + \nonumber \\ & & + B^T\left[\sum_{i=0}^{k-2}(-1)^i(A^T)^i Q A^{k-2-i}\right]B, \quad k \geq 2 .}
We summarize the previous findings in the following theorem:

\begin{theorem}\label{explicitalq}
The constraints $\phi^{(k)}=\sigma^{(k)}\,x+\beta^{(k)}\,p+
\rho^{(k)}\,u$, of the autonomous LQ singular optimal control problem defined by
the matrices $A,~P,\in\R^{n\times n}; ~B,~Q,\in\R^{n\times m};~R\in\R^{m\times m},$
\eqy{\dot{x}&=&Ax+Bu\\L&=&\mitad x^TPx+x^TQu+\mitad u^TRu,}
are given by the following formuli:
\eqy{\label{thbeta} \beta^{(1)} & = &  B^T,  \quad \sigma^{(1)} = -N^T, \quad \rho^{(1)} = -R, \quad \rho^{(2)} = U^1\left(-N^TB+B^TN\right) \\
\beta^{(k)} & = & (-1)^{k}U^{k-1}\cdots U^1B^T(A^T)^{k-1}, \quad k\geq 2, \\
\sigma^{(k)} & = & U^{k-1}\cdots U^1\left(-N^TA^{k-1} + B^T \left[\sum_{i=0}^{k-2}(-1)^i(A^T)^i Q A^{k-2-i}\right]\right),  k\geq 2 ,\\
\rho^{(k)} & = &  U^{k-1} \cdots U^1\left(-N^TA^{k-2}B + (-1)^{k-2}B^T(A^T)^{k-2}N \right.+
\\&  & + \left. B^T \left[\sum_{i=0}^{k-3}(-1)^i(A^T)^i Q A^{k-3-i} \right] B\right), \quad k\geq 3 ,\\
\rho^{(k)} & = & U^{(k)} \Sigma^{(k)} V^{(k)T} ,\\  
\label{thU} U^{k} ~ÊÊ& = & [~0~|~I_{m-r(k)}~]\,U^{(k)\,T} .}
\end{theorem}
If we perturb the matrices $A,B,Q,N,R$ into $A + \delta A$, $B + \delta B$, $Q + \delta Q$, $N + \delta N$ and $R + \delta R$ respectively, then the matrices $\beta^{(k)}$, $\sigma^{(k)}$, $\rho^{(k)}$ will be changed into $\beta^{(k)} + \delta\beta^{(k)}$, $\sigma^{(k)}+ \delta \sigma^{(k)}$, $\rho^{(k)} + \delta \rho^{(k)}$.  Then by using the explicit expressions Eqs. (\ref{thbeta})-(\ref{thU}) it is a straigthforward but tedious computation to obtain the following estimates for the conditioning of the numerical problem of computing the matrices $\beta^{(k)}$, $\sigma^{(k)}$, $\rho^{(k)}$ ($k \geq 2$): 
\eqy{ \frac{|| \delta\widetilde{\beta}^{(k)} ||}{|| \widetilde{\beta}^{(k)}  || }  &\leq& (k-1) \kappa_A \frac{|| \delta A||}{|| A ||} + \kappa_B \frac{|| \delta B ||}{|| B ||}, \\
 \frac{|| \delta\widetilde{\sigma}^{(k)} ||}{|| \widetilde{\sigma}^{(k)}  || }  &\leq& \kappa_A \frac{|| \delta A||}{|| A ||} + \kappa_B \frac{|| \delta B ||}{|| B ||} +  (k-2) \kappa_Q \frac{|| \delta Q ||}{|| Q ||}  + (k-1) \kappa_N  \frac{|| \delta N ||}{|| N ||}, \\
\frac{|| \delta\widetilde{\rho}^{(k)} ||}{|| \widetilde{\rho}^{(k)}  || } &\leq&  C\left[ (k-1) \kappa_A \frac{|| \delta A||}{|| A ||} + \kappa_B \frac{|| \delta B ||}{|| B ||} +  (k-2) \kappa_Q \frac{|| \delta Q ||}{|| Q ||}  + \kappa_N  \frac{|| \delta N ||}{|| N ||} \right] } for some (small) constant $C$.  Notice
that for $k> 1$, the output matrices are not sensitive to variations on the input matrix $R$, which on the other hand is responsible for the launching of the algorithm.  Thus after the first step, the singular matrix singular $R$ dissapears from the computations and does not influence anymore the rest of the construction.

However, in spite of the closed expressions obtained above for the constraints of the system, in order to construct the numerical algorithm to compute them, we will not use Eqs. (\ref{thbeta})-(\ref{thU}) but
rather on we will rely on the recursion Eqs. (\ref{rec1})-(\ref{rec3}). The algorithm will halt whenever at the
step $k$:
\begin{itemize}
\item{$\rho^{(k)}$  is regular, then we can obtain an optimal feedback $u=u(x,p)$ and we
will substitute it in the equations, or,}

\item {$\phi^{(k)}$ is  a linear combination of the previous constraints.  In such a case, there will exists
a ``gauge'' freedom, i.e., some of the controls will be not be determined.}

\end{itemize}
Thus, the scheme of the algorithm will be:

\vspace{0.5cm}
\noindent\begin{minipage}[c]{13cm}
\begin{center}
\textbf{Recursive Constraint Algorithm for singular LQ problems}\vspace{0.3cm}
\hrule
\end{center}\vspace{1mm}
\begin{tabbing}
\texttt{input} $A$, $B$, $Q$, $N$, $R$, $tol$ \\
\emph{Build the constraints matrix:} $\P = [ \sigma^{(1)}~~\beta^{(1)}~~\rho^{(1)}]$\\\\
\texttt{while}  \= $\rank(\rho,tol)$  \emph{is deficient} \& $\rank(\P,~tol)$ \emph{increases}\\
\>$(U,\rho,V)=$ SVD($\rho$); ~\emph{Compute the singular value decomposition of} $\rho^{(k)}$\\
\>\emph{Compute the iterated matrices} $\sigma^{(k+1)},~\beta^{(k+1)},~\rho^{(k+1)}$\\
\>\emph{Build the new constraints matrix:}
$\P=\matrix{ccc}{&\P&\\\sigma^{(k+1)}&\beta^{(k+1)}&\rho^{(k+1)}}$\\
\>\emph{Eliminate the dependent rows of $\P$}\\
\texttt{end while}\\
\texttt{output} $\P,$ $k$\\
\end{tabbing}
\hrule
\end{minipage}

\vspace{0.3cm}
Note that the rank is computed as a numerical rank with tolerance $tol$.
\newpage

\noindent\begin{minipage}[c]{14.5cm}
\begin{center}
\textbf{Pseudocode:  \emph{``final constraint submanifold''}\\ for linear quadratic optimal
control problems}
\vspace{0.3cm}\hrule
\end{center}
\begin{tabbing}
\textbf{input} $A$, $B$, $Q$, $N$, $R$, $tol$\\
$\sigma\leftarrow -N^T;~\beta\leftarrow B^T;~\rho\leftarrow-R;$ ~{Initialize the variables} $\sigma,~\beta,~\rho$\\
$[l,m]\leftarrow$size$(\rho)$; ~{Initialize the dimension $l\times m$ of $\rho$, $l$=maximum value}\\
$\P \leftarrow[\sigma,\beta,\rho];$ ~{Initialize the constraints matrix} $\P$\\
$\P\leftarrow$ independent rows$(\P,~tol)$; ~{Eliminate the dependent rows}\\
$p \leftarrow0$; ~Where $p$ \emph{denotes rank($\P$), it must be 0 to enter in the boucle}\\
$k\leftarrow1;$ \vspace{0.2cm}\\
\\
\textbf{while}  rank\=($\rho,~tol$)$<l~ \& ~$rank($\P,~tol)>$p\\
\>$k\leftarrow k+1;$\\
\>$p\leftarrow$ rank($\P$); ~{Update the rank of $\P$}\\
\>$r\leftarrow$ rank($\rho,~tol$); ~{Update the rank of $\rho^{(k)}$}\\
\>$[l,m]\leftarrow$size($\rho$); ~{Update the dimensions $l\times m$ of $\rho^{(k)}$}\\
\>$[U,\rho,V] \leftarrow $SVD($\rho$); ~{Singular value decomposition of} $\rho^{(k)}$\\
\>$U\leftarrow U^T$\\
\\
\>{Compute the iterated  matrices} $\sigma^{(k+1)},~\beta^{(k+1)},~\rho^{(k+1)}$:\\
\>$\rho\leftarrow~U(r+1:l,:)\cdot[\sigma\,B+\beta\,N]$\\
\>$\sigma\leftarrow~U(r+1:l,:)\cdot[\sigma\,A+\beta\,Q]$\\
\>$\beta\leftarrow ~U(r+1:l,:)\cdot[-\beta\,A^T]$\\
\>$\P\leftarrow\matrix{c}{\P\\\hline \sigma~\beta~\rho}$; ~{Add the new constraints}\\
\>$\P\leftarrow $ independent rows($\P,~tol$); ~{Eliminate the dependent rows}\\
\textbf{end while}\\
\\
\textbf{if} rank$(\P,tol)<=p$\\
    \>$k\leftarrow k-1$\\
\textbf{end if}\\
\\
\textbf{output} ($\P, ~k$)\\
\end{tabbing}\hrule
\end{minipage}

\noindent\begin{minipage}[c]{14.5cm}\vspace{0.3cm}
\begin{center}
\textbf{Pseudocode: \emph{``independent rows''}\\used in the \emph{``final constraint
submanifold''}}
\vspace{0.3cm}\hrule
\end{center}
\begin{tabbing}
\textbf{Procedure} independent rows($\P,~tol$)\\
\\
$[l,c]\leftarrow$ size($\P$); ~Initialize the dimension $l\times c$ of $\P$\\
\\
\textbf{if} $l \geq$ \=$ ~1\&~ \rango (\P,tol)>0$  \textbf{then}\\
\>$F\leftarrow \P(1,:)$\\
\>\textbf{for}  $i=$\= $~2:l$ \\
\>\>\textbf{if} rank\= $(F,tol)< \rango \left( \matrix{c}{F\\\P(i,:)},tol \right)$ \textbf{then}\\
\>\>\>$F\leftarrow\matrix{c}{F\\\P(i,:)}$; ~If the row $i$ is linearly independent we will add it\\
\>\>\textbf{end if}\\
\>\textbf{end for}\\
\textbf{else if}\\
\>$F\leftarrow [\hspace{2mm}]$; ~ If the matrix has zero rank or it is void,  returns the void matrix\\
\textbf{end if}\\
$\P\leftarrow F$\\
\\
\textbf{output} $\P$\\
\end{tabbing}
\hrule
\end{minipage}

\section{Examples and numerical experiments}\label{experiments}

We will discuss here some numerical experiments showing that the numerical algorithm discussed
above behave as expected with respect to stability and consistency.
The microprocessor used for the numerical computations was
Pentium(R), CPU 1.60 GHz, 3.99 MHz, 0.99 GB  RAM, and the program used was MATLAB 7.0.0.

We will describe two types of experiments concerning small ($k = 3$) and large recursive index respectively. 
We are constructing a class of problems that is general enough for the purposes of the numerical stability experiments
we are going to describe and that we can solve and describe the solution explicitly.

In the small index problem $k =3$, we show  that the algorithm is stable with respect to the tolerance used
to compute the numerical rank, $tol$, and with respect to perturbations $\delta$ of the data. We
will also discuss the dependence with the size, $n$, of the matrices.

For the large index ones, we analyze a problem of index  $n-1$, where the algorithm behaves
properly with respect to the number of steps, both regarding the tolerance, $tol$, and the
perturbation of the data, $\delta$.


\vspace{1cm}

{\noindent\bf Small index problems. Small matrices}

Consider the positive semidefinite symmetric $n\times n$ matrix $R$ of
rank 1, thus there will exists an orthogonal matrix $U$ such that  
\eq{\label{svdR} R = U^T R' U} 
such that all elements of $R'$ vanish except $R'_{11} > 0$.  
State and control spaces are both $\R^n$, and the total space $(x,p,u)$ is $\R^{3n}$.
The matrix $A$ is generic and $B$ is an orthonormal matrix, $B^T B = I_n$.  Finally the objective functional
is constructed by using a generic symmetric matrix $Q$, a matrix $N$ of the form $N = BV$ where $V$ is any symmetric matrix and  the matrix $R$ described above.

The primary constraints matrix is given by%
$$ \P^{(1)} = \matrix{c|c|c}{-N^T & B^T & - R}=\matrix{c|c|c}{-I_n & I_n & - R} ,$$%
corresponding to primary constraints%
$$ \phi^{(1)} = -N^T x+ B^Tp  - Ru= 0, $$
Applying the recursive constraint algorithm we obtain:
$$ \dot{\phi}^{(1) }= (-N^TA + B^TQ)x - B^TA^Tp - (B^TN - N^TB)u - RC ,$$
but the SVD of $R$ is given by Eq. (\ref{svdR}), hence the new control coordinate $u^{(1)}$ is given by $u = U^T u^{(1)}$ and
the matrix $U^1$ is just the $(n-1)\times n$ matrix $[0\mid I_{n-1}]U$.  Hence multiplying $\dot{\phi}^{(1)}$ on the left by $U^1$ 
and taking into account that $B^TN - N^TB = 0$, we obtain the
set of secondary constraints:
$$ \phi^{(2)}(x,p,u) =  U^1(-N^TA + B^TQ)x - U^1B^TA^Tp = 0 .$$
Now, computing again the derivative of $\phi^{(2)}$ we obtain the equations:
$$ \dot{\phi}^{(3)} = U^1(-N^TA^2 +  B^TQA - B^TA^TQ)x + U^1B^T(A^T)^2 p + U^1(B^TQB - N^TAB - B^TA^TN )u = 0 .$$ 
We observe that the matrix $B^TQB - N^TAB - B^TA^TN$ will be invertible for generic $A$, $Q$, $B$ and $V$.  For instance
if $A = I$, then $\rho^{(3)}$ reduces to:
$$ \rho^{(3)} = B^TQB - 2V .$$
The algorithm will stop here if $\det (B^TQB - 2V) \neq 0$, this is if $2$ is not an eigenvalue of $V^1B^TQB$.

The numerical experiment  of this problem will consist in applying the algorithm to a
collection of matrices built up as a random perturbation of the matrices $A$, $B$, $Q$, $V$ of size $\delta$,%
$$\tilde{A}=A +\delta A,\quad \|\delta A\|<\delta,\ldots $$ 
It is computed for $n=2,\ldots,202$. We analyze the number of steps before the algorithm stabilizes and compute
the angle, $\alpha$, between the final constraint submanifold of the perturbed problem and the
exact one, this is the error introduced in the problem by the perturbation $\delta A$, etc.
Notice that the original matrix $R$ does not affect higher order constraints, hence its numerical influence 
restricts to launch the algorithm.  Perturbations of $R$ will not affect the computation of higher order constraints
until it will be of the order of $tol$, then the algorithm will stop at the first iteration because if then the system
will be considered to be regular.


\begin{center}
\begin{longtable}{cccccc}
\caption[]{First experiment  (small index, small matrices). $tol=10^{-6}$ }\label{t1}\\
\hline\hline n&$\delta$&\# exact steps&\# steps&codim&$\alpha$/ error\\\hline\hline
\endfirsthead
\caption[]{First experiment  (small index, small matrices). $tol=10^{-6}$}\\
\hline\hline
n&$\delta$&\# exact steps&\# steps&codim&$\alpha$/ error\\\hline\hline
\endhead
\hline\hline
\endfoot
\hline\hline
\endlastfoot
 2&1e-016& 3& 3& 4&0.0000000000000014\\
 2&1e-015& 3& 3& 4&0.0000000000000011\\
 2&1e-014& 3& 3& 4&0.0000000000000059\\
 2&1e-013& 3& 3& 4&0.0000000000000857\\
 2&1e-012& 3& 3& 4&0.0000000000010174\\
 2&1e-011& 3& 3& 4&0.0000000000035392\\
 2&1e-010& 3& 3& 4&0.0000000000213074\\
 2&1e-009& 3& 3& 4&0.0000000002366236\\
 2&1e-008& 3& 3& 4&0.0000000095514555\\
 2&1e-007& 3& 3& 4&0.0000000659481913\\
 2&1e-006& 3& 3& 4&0.0000003955805449\vspace{-1mm}\\
 -- --&-- -- -- -- --&-- --  -- -- -- -- -- --&-- -- -- --&
 -- -- -- --&-- -- -- -- -- -- --  -- -- -- --\vspace{-1mm}\\
 2&1e-005& 3& 1& 2&0.0000048787745951\\
 2&1e-004& 3& 1& 2&0.0000180580284314\\
 2&1e-003& 3& 1& 2&0.0003437294288526\\
 2&1e-002& 3& 1& 2&0.0052691557462037\\
 2&1e-001& 3& 1& 2&0.0558945145125515\\
\hline\hline
n&$\delta$&\# exact steps&\# steps&codim&$\alpha$/ error\\\hline\hline
102&1e-016& 3& 3&304&0.0000000000000071\\
102&1e-015& 3& 3&304&0.0000000000000085\\
102&1e-014& 3& 3&304&0.0000000000000608\\
102&1e-013& 3& 3&304&0.0000000000005422\\
102&1e-012& 3& 3&304&0.0000000000060781\\
102&1e-011& 3& 3&304&0.0000000000579329\\
102&1e-010& 3& 3&304&0.0000000005536516\\
102&1e-009& 3& 3&304&0.0000000060825389\\
102&1e-008& 3& 3&304&0.0000000579413997\vspace{-1mm}\\
 -- --&-- -- -- -- --&-- --  -- -- -- -- -- --&-- -- -- --&
 -- -- -- --&-- -- -- -- -- -- --  -- -- -- --\vspace{-1mm}\\
102&1e-007& 3& 3&302&0.0000004014308009\\
102&1e-006& 3& 3&146&0.0000032681204284\\
102&1e-005& 3& 3&108&0.0000322420378538\vspace{-1mm}\\
 -- --&-- -- -- -- --&-- --  -- -- -- -- -- --&-- -- -- --&
 -- -- -- --&-- -- -- -- -- -- --  -- -- -- --\vspace{-1mm}\\
102&1e-004& 3& 1&102&0.0003320360185045\\
102&1e-003& 3& 1&102&0.0031921082398989\\
102&1e-002& 3& 1&102&0.0340272669549200\\
102&1e-001& 3& 1&102&0.1777792824454151\\
\hline
\hline
n&$\delta$&\# exact steps&\# steps&codim&$\alpha$/ error\\\hline\hline
202&1e-016& 3& 3&604&0.0000000000000078\\
202&1e-015& 3& 3&604&0.0000000000000106\\
202&1e-014& 3& 3&604&0.0000000000000809\\
202&1e-013& 3& 3&604&0.0000000000007860\\
202&1e-012& 3& 3&604&0.0000000000078201\\
202&1e-011& 3& 3&604&0.0000000000801147\\
202&1e-010& 3& 3&604&0.0000000008647247\\
202&1e-009& 3& 3&604&0.0000000081096494\vspace{-1mm}\\
-- --&-- -- -- -- --&-- --  -- -- -- -- -- --&-- -- -- --&
-- -- -- --&-- -- -- -- -- -- --  -- -- -- --\vspace{-1mm}\\
202&1e-008& 3& 3&602&0.0000000577669246\\
202&1e-007& 3& 3&602&0.0000005902470063\\
202&1e-006& 3& 3&262&0.0000048558677195\\
202&1e-005& 3& 3&208&0.0000490952628668\vspace{-1mm}\\
-- --&-- -- -- -- --&-- --  -- -- -- -- -- --&-- -- -- --&
-- -- -- --&-- -- -- -- -- -- --  -- -- -- --\vspace{-1mm}\\
202&1e-004& 3& 1&202&0.0004788619712676\\
202&1e-003& 3& 1&202&0.0046902949060197\\
202&1e-002& 3& 1&202&0.0408456670809488\\
202&1e-001& 3& 1&202&0.2283913263712342\\
\end{longtable}
\end{center}
\vspace{1cm}

\begin{figure}[h!]\begin{center}
\includegraphics[width=10cm]{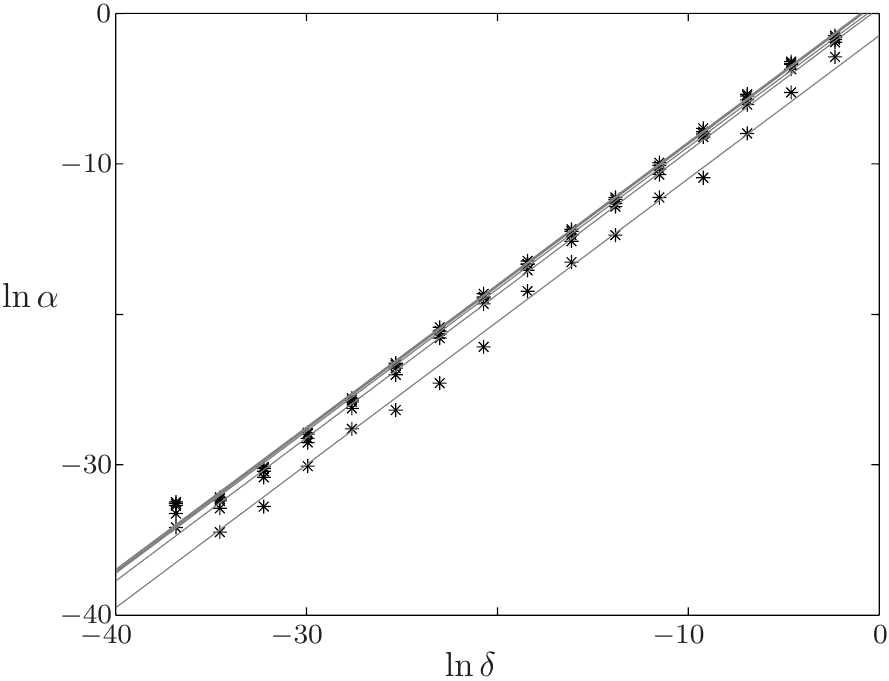}
\caption{Error for the first experiment  (small index, small matrices). $n=2,~52,~102,~152,~202.$
Tolerance used for the computations of the numerical rank equal to $10^{-6}$}
\end{center}\end{figure}

The results show up that the algorithm works well until perturbations of order of
$\delta=10^{-6}$. The least squares approximation of $\ln(\alpha)$
versus $\ln(\delta)$ gives a line of slope  $0.95$, which is consistent with $\alpha=O(\delta).$


In Table \ref{t1} we see that the codimension of the subspace, $codim=3n-1$, that is, the number of
rows  of the constraints matrix fails at $n=2$ when $\delta$ is of the order of $tol$. However,
when $n$ grows $codim$ fails for smaller $\delta$.

Moreover the results show that the algorithm is insensitive to the size of the original matrices.
In fact, if we select a fixed value of the perturbation, $\delta=10^{-6}$, and analyze the error
for different values of $n$, we obtain   Table \ref{t2}.

\begin{center}
\begin{table}
\caption{First experiment  (small index, small matrices). $\delta=10^{-9},$
$tol=10^{-6}$}\label{t2}
\begin{tabular}{cccccc}\hline\hline
n&$\delta$&\# exact steps&\# steps&codim&$\alpha/$error\\\hline\hline
2&1e-009& 3& 3& 4&0.0000000009585285\\
22&1e-009& 3& 3&64&0.0000000028931919\\
42&1e-009& 3& 3&124&0.0000000034146838\\
62&1e-009& 3& 3&184&0.0000000046018755\\
82&1e-009& 3& 3&244&0.0000000053815907\\
102&1e-009& 3& 3&304&0.0000000055767568\\
122&1e-009& 3& 3&364&0.0000000064274692\\
142&1e-009& 3& 3&424&0.0000000070707452\\
162&1e-009& 3& 3&484&0.0000000072378230\\
182&1e-009& 3& 3&544&0.0000000077680372\\
202&1e-009& 3& 3&604&0.0000000083969815\\
\hline\hline
\end{tabular}
\end{table}
\end{center}
Again data obtained indicates  heuristically that $\alpha=O(\sqrt{n})$. The least squares
approximation of $\ln(\alpha)$ versus $\ln(n)$ gives us an slope of  $0.47$.

\vspace{1cm}

\begin{center}\begin{figure}[h!]
\includegraphics[width=10cm]{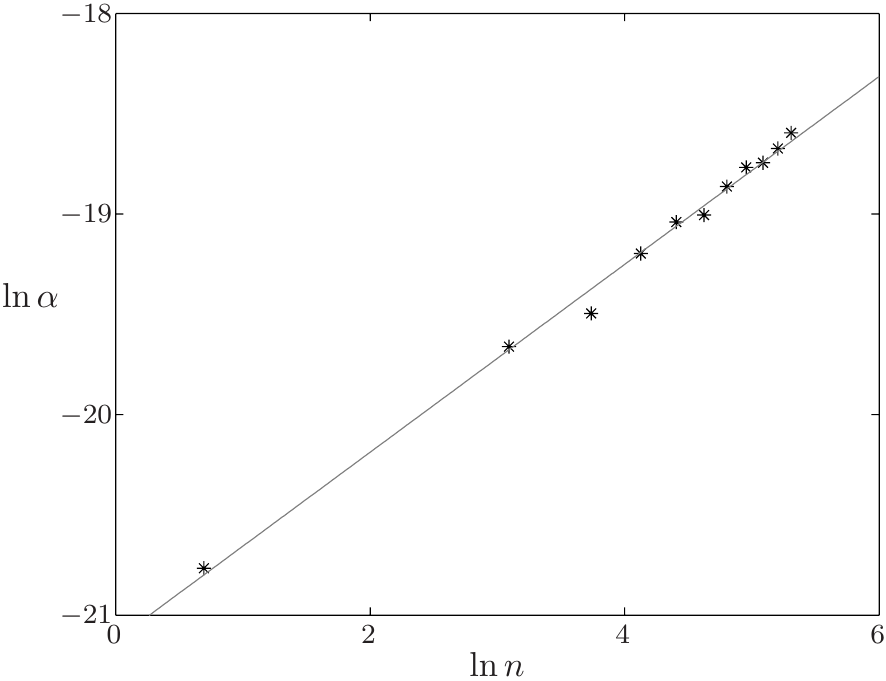}
\caption{Error for the first experiment  (small index, small matrices).  $n=2\to~202$,
$\delta=10^{-6}$, $tol=10^{-6}$}
\end{figure}\end{center}

\vspace{1cm}

\newpage

{\noindent\bf Small index problems. Large matrices}

Let us consider the linear--quadratic problem where $A$ is proportional to the identity matrix,
$A=\alpha I$, where $\alpha\in\R$.

The primary constraint will be
\eq{\phi^{(1)}=\sigma^{(1)}x+\beta^{(1)}p+\rho^{(1)}u=-N^Tx+B^Tp-Ru,}the general form of the
$(k+1)-th$ constraint will be \eqi{\phi^{(k+1)}&=&\sigma^{(k+1)}x+\beta^{(k+1)}p+\rho^{(k+1)}u=
[~0~~|~I_{m-r_{k}}~]\,U^{(k)\,T}\left(\sigma^{(k)}\dot{x}+\beta^{(k)}\dot{p}\right)=\\
&=&[~0~|~I_{m-r_k}~]\,U^{(k)\,T}\left[(\alpha\sigma^{(k)}+ \beta^{(k)}Q)x+(-\alpha\beta^{(k)})p+
(\sigma^{(k)}B+\beta^{(k)}N)u\right].}%
Here we use the same notation as in Section \ref{pcasect}, that is, after applying the SVD to the matrix
$\sigma^{(k)}$, we obtain $\rho^{(k)}=U^{(k)}\,\Sigma^{(k)}\,V^{(k)\,T}$, and
$U_k:=[~0~~|~I_{m-r_k}~]\,U^{(k)\,T}$. Computing the constraints, we obtain
\eqi{\phi^{(2)}&=&\sigma^{(2)}x+\beta^{(2)}p+\rho^{(2)}u=
U_1\left[(\alpha\sigma^{(1)}+ \beta^{(1)}Q)x+(-\alpha\beta^{(1)})p+
(\sigma^{(1)}B+\beta^{(1)}N)u\right]=\\
&=&U_1\left[(-\alpha N^T+ B^T Q)x+(-\alpha B^T)p+(-N^TB+B^T N)u\right],\\
\phi^{(3)}&=&\sigma^{(3)}x+\beta^{(3)}p+\rho^{(3)}u=
U_2\left[(\alpha\sigma^{(2)}+ \beta^{(2)}Q)x+(-\alpha\beta^{(2)})p+
(\sigma^{(2)}B+\beta^{(2)}N)u\right]=\\
&=&U_2U_1\left[(-\alpha^2 N^T+ \alpha B^T Q-\alpha B^T Q)x+(\alpha^2 B^T)p+(-\alpha N^TB+ B^T QB-\alpha B^TN)u\right]
=\\
&=&U_2U_1\left[(-\alpha^2 N^T)x+(\alpha^2 B^T)p+(-\alpha (N^TB+B^TN)+ B^T QB)u\right],\\
\phi^{(4)}&=&\sigma^{(4)}x+\beta^{(4)}p+\rho^{(4)}u=
U_3\left[(\alpha\sigma^{(3)}+ \beta^{(3)}Q)x+(-\alpha\beta^{(3)})p+
(\sigma^{(3)}B+\beta^{(3)}N)u\right]=\\
&=&U_3U_2U_1\left[(-\alpha^3 N^T+ \alpha^2 B^TQ)x+(-\alpha^3 B^T)p+
(-\alpha^2 N^TB+\alpha^2 B^TN)u\right].}
So the constraints matrix will look as
\eqi{\footnotesize\P=\matrix{c|c|c}{-N^T&B^T&-R\\ &&\\
U_1\left[-\alpha N^T+B^TQ\right]&-U_1\left[\alpha B^T\right]&U_1\left[-N^TB+B^TN\right]\\ &&\\
-U_2U_1\left[\alpha^2N^T\right]&U_2U_1\left[\alpha^2B^T\right]&
U_2U_1\left[\alpha(-N^TB-B^TN)+B^TQB\right]\\ &&\\
U_3U_2U_1\alpha^2\left[-\alpha N^T+B^TQ\right]&-U_3U_2U_1\alpha^2\left[\alpha B^T\right]&
U_3U_2U_1\alpha^2\left[-N^TB+B^TN\right]}.}

We can see that the fourth row is related with the second one by
$\mathrm{row}_4=U_3U_2\alpha^2\mathrm{row}_2$, so the algorithm will stop here if it did not do it before.

For the numerical implementation we choose the following matrices: $Q=A=I_n \in\R^{n\times n},$
$B^T=(1,\ldots,1)\in\R^{n\times 1},~N^T=(0,\ldots,0)\in\R^{n\times 1},~R=0;$ so the constraints
matrix will have only three rows
\eqi{\P=\matrix{c|r|c}{0,\ldots,0&1,\ldots,\hspace{2.6mm}1&0\\1,\ldots,1&-1,\ldots,-1&0\\
0,\ldots,0&1,\ldots,\hspace{2.6mm}1&n},}where $n$ is the dimension of the matrices $A$ and $Q$.
Thus in the third row we obtain  optimal feedback and the final constraint submanifold is given by
the following equations: $x_1+\cdots +x_n=p_1+\cdots+ p_n=u=0$.

We apply the numerical algorithm  for the previows matrices for $n=1000$, tolerance equal to
$10^{-16}$ and we compare the solution obtained with the perturbed  matrices: $\tilde{A}=A+\delta
A,$ $\|\delta A\|<\delta,$ $\tilde{N}=N+\delta N,$ $\|\delta N\|<\delta$ and $\tilde{B}=B+\delta
B,$ $\|\delta B\|<\delta$, where $\delta=10^{-16}\to 10^{-1}$. Again, as the final constraint
submanifold of the original problem and the perturbed one must be the same, we measure the angle
between them, this is going to be the error, and we show it in Table \ref{t3}.

\begin{center}
\begin{table}\caption{Second experiment (small index, large matrices).  $tol=10^{-16}$}\label{t3}
\begin{tabular}{cccccc}\hline\hline
n&$\delta$&\# exact steps&\# steps&codim&$\alpha$/ error\\\hline\hline
1000&1e-016& 3& 3& 3&0.00000000000002\\
1000&1e-015& 3& 3& 3&0.00000000000005\\
1000&1e-014& 3& 3& 3&0.00000000000010\\
1000&1e-013& 3& 3& 3&0.00000000000088\\
1000&1e-012& 3& 3& 3&0.00000000000921\\
1000&1e-011& 3& 3& 3&0.00000000008945\\
1000&1e-010& 3& 3& 3&0.00000000094362\\
1000&1e-009& 3& 3& 3&0.00000000895734\\
1000&1e-008& 3& 3& 3&0.00000009224712\\
1000&1e-007& 3& 3& 3&0.00000089770415\\
1000&1e-006& 3& 3& 3&0.00000934812796\\
1000&1e-005& 3& 3& 3&0.00008953156206\\
1000&1e-004& 3& 3& 3&0.00088056539460\\
1000&1e-003& 3& 3& 3&0.00953194630784\\
1000&1e-002& 3& 3& 3&0.10128762315235\\
1000&1e-001& 3& 3& 3&0.80931668976548\\
\hline\hline
\end{tabular}
\end{table}
\end{center}

Again the data shows that $\alpha=O(\delta)$ and, consistently with the previous results, the slope
of the least squares approximation of $\ln(\alpha)$ versus $\ln(\delta)$ is $0.96$.

\newpage

\begin{figure}[h!]\begin{center}
\includegraphics[width=10cm]{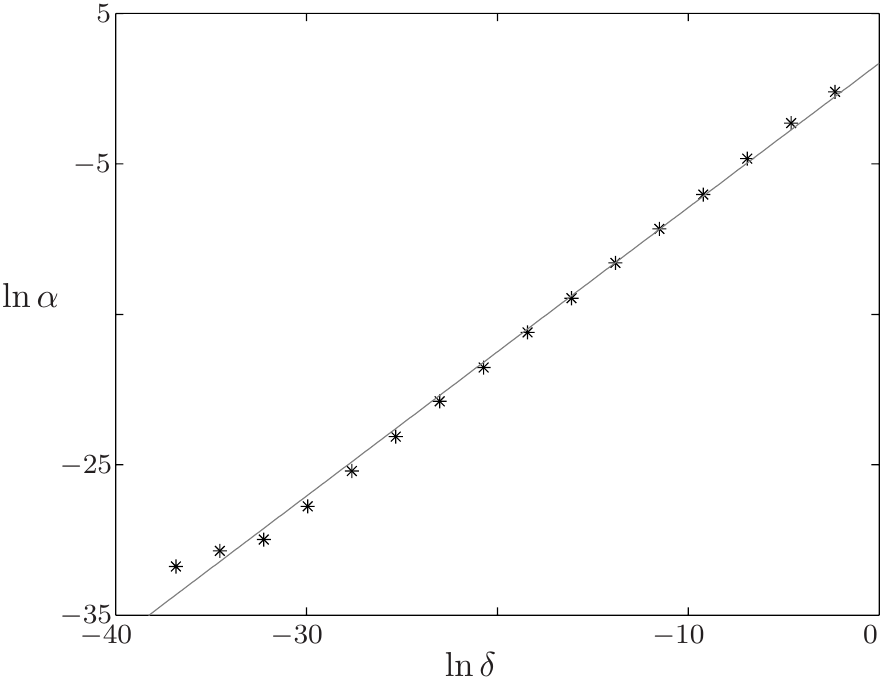}
\caption{Error for the second experiment (small index, large matrices).
$n=1000$.  $tol=10^{-16}$}
\end{center}\end{figure}

{\noindent\bf Large index problems}

Consider the following problem:%
\eq{ A\in\R^{n\times n}, ~Q=A+A^T, ~B\in\R^{n\times 1}, ~N=B, R=0.}
Computing the matrices $\rho^{(1)},~\rho^{(2)},\ldots,~\rho^{(k)}$, we get
\eqi{\rho^{(1)} & = & R=0, \\
     \rho^{(2)} & = & B^TN-N^TB = B^TB-B^TB = 0, \\
     \rho^{(3)} & = & -N^TAB - B^TA^TN + B^TQB = B^T[-A - A^T + A + A^T]B = 0, \\
     \rho^{(4)} & = & B^T[-A^T + (-A^T)^2 + (A + A^T)A - A^T(A + A^T)]B = 0, \\
                & \vdots & \\
    \rho^{(k+1)}& = & B^T[-A^{k-1} + (-1)^{k-1}(A^T)^{k-1} + \sum_{i=0}^{k-2}(-1)^i(A^T)^i(A+A^T)A^{k-2-i}]B= \\
                & = & B^T[-A^{k-1} + (-1)^{k-1}(A^T)^{k-1} \\
                & + & \sum_{j=1}^{k-1}-(-1)^j(A^T)^jA^{k-1-j} + \sum_{j=0}^{k-2}(-1)^j(A^T)^jA^{k-1-j}]B= \\
                & = & B^T[-A^{k-1} + (-1)^{k-1}(A^T)^{k-1} - (-1)^{k-1}(A^T)^{k-1} + A^{k-1}]B = 0.}%
We obtain that these matrices are always zero and there will not exist optimal feedback. Let us
consider now the remaining matrices
\eqi{ \sigma^{(1)} & = & -B^T ,\\
      \sigma^{(2)} & = & -B^TA + B^T(A + A^T) = B^TA^T, \\
      \sigma^{(3)} & = & B^T[-A^2 + (A + A^T)A - A^T(A+A^T)] = -B^T(A^T)^2,\\
                   &\vdots&\\
     \sigma^{(k+1)}& = & B^T[-A^k+\sum_{i=0}^{k-1}(-1)^i(A^T)^i(A+A^T)A^{k-1-i}]=\\
                   & = & B^T[-A^k+\sum_{i=0}^{k-1}(-1)^i(A^T)^{i+1}A^{k-1-i}+\sum_{i=0}^{k-1}(-1)^i(A^T)^{i}A^{k-i}]=\\
                   & = & B^T[-A^k+\sum_{j=1}^{k}-(-1)^j(A^T)^{j}A^{k-j}+\sum_{j=0}^{k-1}(-1)^j(A^T)^{j}A^{k-j}]=\\
                   & = & B^T[-A^k-(-1)^k(A^T)^{k}+A^k=(-1)^{k-1}B^T(A^T)^{k}.}
\eqi{\beta^{(1)}& = & B^T,\\
     \beta^{(2)}& = & -B^TA^T,\\
                &\vdots&\\
   \beta^{(k+1)}& = & (-1)^kB^T(A^T)^k.}%
Hence, we will obtain
\eqy{\beta^{k+1}=-\sigma^{k+1}=(-1)^kB^T(A^T)^k,~~~\rho^{k+1}=0.}%
and the constraints matrix will look as follows
\eqi{\P &=& \matrix{cccc}{-B^T & B^T & & 0\\ B^TA^T & -B^TA^T & & 0 \\-B^T(A^T)^2 & B^T(A^T)^2 & &0\\
\vdots & \vdots & & \vdots\\ -(-1)^kB^T(A^T)^k & (-1)^kB^T(A^T)^k & & 0}.}%
Moreover, it is clear from the previous considerations that the algorithm will stop only when at a
given step we will obtain a linear combination of the previous rows.   Suppose that the minimum
polynomial of the matrix $A$ is of degree $q$, then there are two possibilities for the algorithm
to stop:

\begin{itemize}

\item $B\notin~\ker(A^l),~~l=1,\ldots,q$, then the row $q+1$ is a linear combination
of the previous ones.

\item $B\in~\ker(A^p),~~0<p\leq q$, then the row $p+1$ vanish.

\end{itemize}


We apply the algorithm to a problem where the pair $(A,B)$ is such that $A\in\R^{n\times n}$ is a
nilpotent matrix of index $n$, i.e.,  $A^{n-1}\neq 0$, $A^n=0$, and $B\notin~\ker(A^l)$,
$l=1,\ldots,n-1$, i.e., $B^T=(1,\ldots,1)\in\R^{1\times n}$.  The index of the algorithm is $k =
n$. Perturbing the matrices as: $\tilde{A}=A+\delta A,$ $\|\delta A\|<\delta,$ $ \tilde{Q} =
\tilde{A}+ \tilde{A}^T,$ $\tilde{B}=B+\delta B,$ $\|\delta B\|<\delta$ , $\tilde{N}=\tilde{B}$ and
$\tilde{R}=\delta R,$ $\|\delta R\|<\delta;$
with $\delta=10^{-16}\to 10^{-5}$, we get for $n=20$  with tolerance equal to $10^{-6}$  Table
\ref{t4}.


\begin{center}
\begin{table}[h!]\caption{Third experiment (large index).  $tol=10^{-6}$}
\label{t4}
\begin{tabular}{cccccc}\hline\hline
n&$\delta$&\# exact steps&\# steps&codim&$\alpha$/ error\\\hline\hline
20&1e-016&20&20&20&0.00000000000001\\
20&1e-015&20&20&20&0.00000000000002\\
20&1e-014&20&20&20&0.00000000000001\\
20&1e-013&20&20&20&0.00000000000004\\
20&1e-012&20&20&20&0.00000000000021\\
20&1e-011&20&20&20&0.00000000000652\\
20&1e-010&20&20&20&0.00000000001940\\
20&1e-009&20&20&20&0.00000000034574\\
20&1e-008&20&20&20&0.00000000647005\\
20&1e-007&20&20&20&0.00000006884601\\
20&1e-006&20&20&20&0.00000047329442\\
20&1e-005&20& 1& 1&0.00000157131430\\
\hline\hline
\end{tabular}
\end{table}
\end{center}

In this experiment, the value of the slope of the least squares approximation of the data is $0.84$.

\newpage

\begin{center}
\begin{figure}[h!]
\includegraphics[width=10cm]{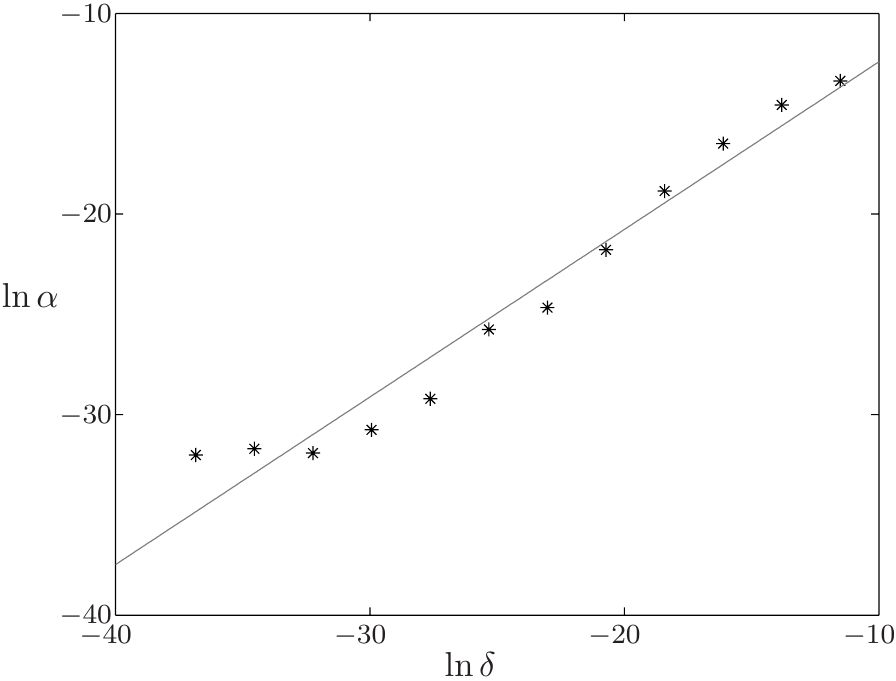}
\caption{Error for the third experiment (large index). $n=20$,  $tol=10^{-6}$}
\end{figure}
\end{center}



\end{document}